\documentclass[a4in,12pt]{article}
\usepackage{amsmath}
\usepackage{amssymb}
\usepackage{amscd}

\newtheorem{theoreme}{Th\'eor\`eme}[section]
\newtheorem{proposition}[theoreme]{Proposition}

\newtheorem{lemme}[theoreme]{Lemme}

\newtheorem{exemple}[theoreme]{Exemple}
\newtheorem{remarque}[theoreme]{Remarque}

\newenvironment{preuve}{\begin{trivlist} \item[]{\it Preuve---}}
{\par\hfill $\square$\end{trivlist}}
\renewcommand{\P}{\mathbb{P}}
\newcommand{\C}{\mathbb{C}}

\newcommand{\R}{\mathbb{R}}
\newcommand{\N}{\mathbb{N}}

\newcommand{\K}{{\cal K}}

\newcommand{\ddc}{{\rm dd^c}}
\renewcommand{\d}{{\rm d}}

\newcommand{\cad}{{\it c.-\`a-d. }}
\newcommand{\ie}{{\it i.e. }}

\newcommand{\Lone}{{{\rm L}^1}}

\newcommand{\Linfty}{{{\rm L}^\infty}}

\newcommand{\Loneloc}{{{\rm L}^1_{\rm loc}}}

\newcommand{\supp}{{\rm supp}}
\newcommand{\E}{{\cal E}}
\newcommand{\D}{{\cal D}}

\newcommand{\voir}{{\it voir }}

\renewcommand{\o}{{\rm o}}

\title{Dynamique des applications polynomiales semi-r\'eguli\`eres}
\author{Tien-Cuong Dinh et Nessim Sibony}
\date{Nouvelle version}
\begin{document}
\maketitle
\begin{abstract} 
For any proper polynomial map $f:\C^k\longrightarrow \C^k$ define the
function $\alpha$ as 
$$\alpha(z):=\limsup_{n\rightarrow\infty} 
\frac{\log^+\log^+|f^n(z)|}{n} \mbox{  where  } \
\log^+:=\max(\log, 0).$$ 
Let $f=(P_1,\ldots,P_k)$ be a proper polynomial map. 
We define a notion of $s$-regularity using the extension of $f$ to
$\P^k$. When $f$ is (maximally) regular
we show that 
the function $\alpha$ is l.s.c and takes only finitely many values:
  $0$ and $d_1$, $\ldots$, $d_k$, where $d_i:=\deg P_i$. 
We then describe dynamically the sets $\{\alpha\leq
    d_i\}$. If $d_i>1$, 
this allows us to construct the equilibrium measure
    $\mu$ associated to $f$ as a generalized intersection of positive
    currents. We then gives an estimate of the Hausdorff dimension of
    $\mu$. This is a special case of our results. We extend the
    approach to the larger class of $(\pi,s)$-regular
    maps. This gives an understanding of the biggest values of
    $\alpha$. The results can be applied to construct dynamically
    interesting measures for automorphisms.
\end{abstract}
\section{Introduction}
Soit $f=(P_1,\ldots,P_k)$ une application polynomiale propre de $\C^k$
dans $\C^k$. Quitte \`a conjuguer $f$ par une permutation des
coordonn\'ees on peut supposer que $\deg P_1\geq \cdots \geq\deg
P_k$. D\'efinissons des entiers $l_i$ v\'erifiant $1=l_0<l_1<\cdots
<l_m=k+1$ tels que les composantes de
$$f_{(i)}:=P_{(i)}=(P_{l_{i-1}},\ldots, P_{l_i-1})$$
soient de m\^eme degr\'e $d_i$ avec $d_1>d_2>\cdots >d_m$. 
Notons $P_{(i)}^+$ la partie homog\`ene de plus haut degr\'e de
$P_{(i)}$. 
Notons
\'egalement $f$ l'extension de $f$ comme application m\'eromorphe \`a
$\P^k$ dont $[z_1:\cdots:z_k:t]$ sont les coordonn\'ees homog\`enes.
\par
Pour tout $z\in\C^k$ d\'efinissons la constante $\alpha(z)\geq 0$ par 
$$\alpha(z):=\limsup_{n\rightarrow\infty} 
\frac{\log^+\log^+|f^n(z)|}{n}$$ 
o\`u $\log^+:=\max(\log, 0)$.
\par
En g\'en\'eral, $\alpha$ prend une infinit\'e de valeurs. Lorsque
$\alpha(z)>1$, $\alpha(z)$ repr\'esente la vitesse d'\'echapement de
$f^n(z)$ vers l'infini. Lorsque $f$ se prolonge holomorphiquement \`a
l'infini dans $\P^k$, auquel cas $m=1$ et $\{P^+_{(1)}=0\}$ est
r\'eduit \`a l'origine, la fonction $\alpha$ ne prend que deux valeurs
$0$ et $d_1$. 
Pour les applications r\'eguli\`eres que nous introduisons,
nous montrons que $\alpha$ est s.c.i. et ne prend
que les valeur $d_1$, $\ldots$, $d_m$, $0$.
\par
On introduit des fonctions de Green partielles sur les ferm\'es
$\K_i:=\{\alpha\leq d_i\}$ par 
$$G_i(z):=\lim_{n\rightarrow 0} \frac{\log^+|f^n(z)|}{d_i^n}.$$
On montre que si $i<m$ ou si $i=m$ et $d_m>1$,
la fonction $G_i$ est continue sur $\K_i$. Lorsque
$d_m>1$, on obtient la mesure d'\'equilibre $\mu$ comme produit
d'intersection g\'en\'eralis\'e de courants positifs d\'efinis \`a
l'aide des fonctions $G_i$. L'\'etude des fonctions $G_i$ permet
d'obtenir une estimation de la dimension de Hausdorff de $\mu$.
\par
Lorsque $f$ se prolonge holomorphiquement \`a l'infini et $d_1>1$ on a
$m=1$. La fonction $G=\lim
d_1^{-n}\log^+|f^n|$ est d\'efinie partout et $\mu=(\ddc G)^k$. 
\par
Nous renvoyons en particulier \`a \cite{BriendDuval2, DinhSibony2,
  FornaessSibony2, Sibony2, HubbardPapadopol} pour divers aspects de la
  dynamique de ces applications. Lorsque $f$ est r\'eguli\`ere et
  $d_m>1$, $f$ est \`a allure polynomiale au sens de
  \cite{DinhSibony2} et la mesure $\mu$ que nous construisons ici est
  la m\^eme que dans \cite{DinhSibony2} (\voir le paragraphe 4).
\par
Nous nous int\'eressons dans cet article \`a des classes plus
g\'en\'erales que les applications r\'eguli\`eres \`a savoir les
applications $s$-r\'eguli\`eres et $(\pi,s)$-r\'eguli\`eres. Pour
d\'ecrire ces classes introduisons quelques notations.
\par
$$z_{(i)}:=(z_{l_{i-1}},\ldots, z_{l_i-1}) \ \ \ \  
z^d_{(i)}:=(z^d_{l_{i-1}},\ldots, z^d_{l_i-1})$$
$$z_{(<i)}:=(z_{(1)},\ldots, z_{(i-1)}) \ \ \ \ 
z_{(\leq i)}:=(z_{(1)},\ldots, z_{(i)})$$
$$z_{(>i)}:=(z_{(i+1)},\ldots, z_{(m)}) \ \ \ \ 
z_{(\geq i)}:=(z_{(i)},\ldots, z_{(m)})$$
$$|z|_{(i)}:=|z_{(i)}| \ \ \ \mbox{ etc.}$$
On dira que $g$ et $h$ sont 
{\it comparables} 
quand $z\rightarrow X$ et on
notera $g(z)\sim h(z)$ s'il existe des
constantes $0<c<c'$ telles que $cg(z)\leq h(z)\leq c' g(z)$ pour $z$
suffisamment proche de $X$. 
On dira que $g$ et $h$ sont {\it \'equivalents}
quand $z\rightarrow X$ et on
note $g(z)\simeq h(z)$ si $g(z)/h(z)$ tend vers 1
quand $z\rightarrow X$. 
\par
On identifie $\C^k$ \`a $\P^k\setminus\{t=0\}$ et pour toute
application polynomiale $Q:\C^k\longrightarrow \C^n$, $\{Q=0\}$
d\'esigne le sous ensemble alg\'ebrique de $\P^k$ adh\'erence de
$Q^{-1}(0)$.
\par
Posons $I_0:=\{t=0\}$ et $X_0:=\emptyset$. 
Pour tout $1\leq i\leq m$, posons
$$I_i:=I_{i-1}\cap \big\{P^+_{(i)}=0\big\}$$
et
$$X_i:=I_{i-1}\cap \big\{z_{(> i)}=0\big\}.$$
On dira que $f$ est {\it $s$-r\'egulier}
si $I_i\cap X_i=\emptyset$
pour tout $1\leq i\leq s$ et que $f$ est {\it r\'egulier}
s'il
est $m$-r\'egulier. Si $f$ est $s$-r\'egulier, d'apr\`es le
th\'eor\`eme de B\'ezout, on a n\'ecessairement 
$\dim X_i=l_i-l_{i-1}-1$ et $\dim I_i=k-l_i$ pour tout $1\leq i\leq
s$. L'ensemble $X_i$ appara\^{\i}t comme l'image de $I_{i-1}\setminus I_i$
par une restriction convenable de $f$ \`a $I_{i-1}$. L'ensemble $I_i$
appara\^{\i}t comme l'ensemble d'ind\'etermination de cette restriction de
$f$ \`a $I_{i-1}$.
\par
Si $f$ est $1$-r\'egulier, on a
$I_1\cap X_1=\emptyset$. En particulier,
$f$ est {\it alg\'ebriquement stable}
\cite{FornaessSibony2,Sibony2}, \cad qu'aucune
hypersurface n'est envoy\'ee par un it\'er\'e de $f$, dans $I_1$. Il en
r\'esulte que
le degr\'e alg\'ebrique
de $f^n$ est \'egal \`a $d_1^n$. On peut alors d\'efinir la
{\it fonction de Green}
$G_1$ par
$$G_1(z):=\lim_{n\rightarrow \infty}\frac{\log^+|f^n(z)|}{d_1^n}.$$
C'est une fonction p.s.h. sur $\C^k$. Elle d\'efinit un courant
positif ferm\'e $T_1:=\ddc G_1$ qui se prolonge \`a $\P^k$
car $G_1(z)-\log^+|z|$ est born\'ee sup\'erieurement. 
D'apr\`es la proposition 5.3,
ce courant ne charge pas les ensembles pluripolaires
de $\P^k\setminus I_1$ car il admet localement un potentiel
born\'e en tout point de $\P^k\setminus I_1$. Il ne peut non plus  
charger $I_1$ car $I_1$ est de codimension au moins $2$ dans $\P^k$. 
\par
De fa\c con g\'en\'erale, on pose
$G_0:=0$, $\K_0=\C^k$ et
lorsque les expressions suivantes ont un sens,
on pose 
$$G_{i,n}(z):=\frac{\log^+|f^n(z)|}{d_i^n}$$
$$G_i(z):=\lim_{n\rightarrow\infty} G_{i,n}(z).$$
On d\'efinit 
$$U_i:=\big\{z\in\C^k,\  f^n(z) \mbox{ tend vers } X_i \big\}$$
et
$$\K_i:=\K_{i-1}\setminus U_i$$
pour tout $1\leq i\leq m$. 
%
%
\par
Si $f$ est $s$-r\'egulier et $1\leq i\leq s$, 
nous montrons que pour $z\in U_i$, $\alpha(z)=d_i$ et on pr\'ecise la
dynamique de $f$ dans $U_i$ et son compl\'ementaire. Pr\'ecisons cela.
\par
Si $f$ est
1-r\'egulier, $U_1$ est le bassin de $X_1$. 
On v\'erifiera que $U_1$, $\K_1$ sont
invariants par $f$ et $f^{-1}$,
$\overline\K_1\subset \K_1\cup I_1$, $X_1\cap
I_1=\emptyset$ et $\dim I_1=k-l_1$ o\`u $I_1$ est l'ensemble
d'ind\'etermination de $f$.
La fonction de Green $G_1$ pr\'ecise l'\'echapement vers l'infini.
On montrera que $G_1$ est continue, positive, nulle exactement sur
$\K_1$ et \`a croissance logarithmique \`a l'infini.
Elle est de plus invariante par $f$: $G_1\circ f=d_1 G_1$.
Pour tout
$1\leq j\leq l_1-1$, le courant $T_j:=(\ddc G_1)^j$ est positif,
ferm\'e, de bidegr\'e $(j,j)$,
invariant par $f$ et ne charge pas
les ensembles pluripolaires. Le courant $T_{l_1-1}$ est port\'e
par $\overline \K_1$.
\par
Il s'agit maintenant d'analyser la dynamique de la restriction de
$f$ \`a $\K_1$. 
Supposons que $f$ est 2-r\'egulier. Le sous ensemble analytique $X_2$
de  $I_1$ est attirant; il est de dimension
$l_2-l_1-1$. Si $z$ appartient \`a un petit voisinage $V_2$
de $X_2$ dans $\K_1\cup I_1$, on a 
$$c^{-1}|z|^{d_2}\leq |f(z)|\leq c|z|^{d_2}$$
o\`u $c>0$ est une constante. Le bassin $U_2$
est \'egal \`a $\bigcup_{n\geq 0}f^{-n}(V_2)$; son compl\'ementaire
$\K_2$ v\'erifie $\overline \K_2\subset \K_2\cup I_2$. Rappelons que
$I_2\subset I_1$
est un sous-ensemble analytique de dimension $k-l_2$ v\'erifiant
$X_2\cap I_2=\emptyset$. La deuxi\`eme fonction de Green
$G_2(z)$
est finie et continue, positive sur $\K_1$, elle est 
\'egale \`a $+\infty$ sur $U_1$,
nulle exactement sur
$\K_2$. Sur $\K_1$, elle a une croissance
logarithmique \`a l'infini.
On a la relation invariante $G_2\circ f=d_2 G_2$.
Pour tout
$l_1\leq j\leq l_2-1$, le courant $T_j:=(\ddc G_2)^{j-l_1+1}\wedge T_{l_1-1}$
est positif,
ferm\'e, de bidegr\'e $(j,j)$, 
invariant par $f$ et ne charge pas
les ensembles pluripolaires. Le courant $T_{l_2-1}$ est port\'e
par $\overline \K_2$.
\par
Suivant l'ordre $s$ de la r\'egularit\'e, la construction peut se
poursuivre. Lorsque $s=m$, on trouver les applications r\'eguli\`eres.
Au paragraphe 2, nous explicitons le cas des applications
$s$-r\'eguli\`eres. Au paragraphe 3, 
nous \'etendons la th\'eorie aux applications $(\pi,s)$-r\'eguli\`eres
et nous donnons une estimation de la dimension de Hausdorff de
$\mu$. 
Au paragraphe 4, nous r\'esumons d'autres propri\'et\'es dynamiques
des applications r\'eguli\`eres et $\pi$-r\'eguli\`eres. Observons que
lorsqu'on fixe $d_1>d_2>\cdots>d_m$, dans l'espace de param\`etres,
les familles d'applications r\'eguli\`eres et semi-r\'eguli\`eres
sont des ouverts Zariski denses.
Nous avons rassembl\'e dans un appendice les propri\'et\'es des
fonctions p.s.h. par rapport \`a un courant positif ferm\'e, que nous
utilisons. 
\par
Dans \cite{GuedjSibony}, on trouve d\'ej\`a la d\'efinition de
fonctions de Green partielles pour certaines automorphismes de $\C^k$ et
pour des endomorphismes de $\C^2$ dans \cite{FavreGuedj}. On trouve
\'egalement dans \cite{GuedjSibony,Guedj} une notion de faible
r\'egularit\'e voisine de la $1$-r\'egularit\'e. L'int\'er\^et de
notre approche ici est que nous d\'eduisons les estimations
n\'ecessaires \`a la construction des fonctions de Green partielles 
d'hypoth\`eses g\'eom\'etriques faciles \`a
v\'erifier. Nous donnons en particulier (proposition 3.1) une
caract\'erisations des applications $f$ $(\pi,s)$-r\'eguli\`eres de
$\C^2$ en termes des polyg\^ones de Newton des composantes de
$f$. Notre construction permet d'obtenir des mesures invariantes
int\'eressantes dans le cas des
automorphismes polynomiaux (remarque 3.9). 
\section{Endomorphismes r\'eguliers}
Dans la suite, notons $d_t$ le degr\'e topologique de $f$, \cad le
nombre de pr\'eimages d'un point $z\in\C^k$ compt\'ees avec
multiplicit\'e. Le degr\'e topologique ne d\'epend pas du point
$z$. Notons $\K$ l'ensemble des points d'orbite born\'ee. 
{\it L'exposant de Lojasiewicz} de $f^n$ sera not\'e $\lambda_n$. C'est
la meilleure constante positive v\'erifiant $|f^n(z)|\geq
c|z|^{\lambda_n}$ pour $z$ assez grand o\`u $c>0$ est une constante. 
Cette constante $\lambda_n$ existe toujours \cite{Ploski}. La suite
$(\lambda_n^{1/n})$ d\'ecroit vers une constante $\lambda_\infty$ que
nous appelons
{\it l'exposant de Lojasiewicz asymptotique} de $f$.
On a le th\'eor\`eme suivant:  
\begin{theoreme} Soit $f:\C^k\longrightarrow \C^k$ un
endomorphisme polynomial propre $s$-r\'egulier. On suppose que
$1\leq s\leq m-1$ ou bien que $s=m$ et $d_m\geq 2$.
Alors
pour tout $1\leq i\leq s$, il existe une suite de nombres r\'eels
positifs 
$(c_{i,n})\rightarrow 0$ telle que la suite de fonctions
$G_{i,n}+c_{i,n}$
d\'ecroit sur $\K_{i-1}$ vers une fonction $G_i$ continue,
invariante: $G_i\circ f=d_i G_i$. De plus, la fonction $G_i(z)-\log^+|z|$
est continue sur $\overline\K_{i-1}\setminus I_i$ et
\begin{eqnarray*}
\K_{i-1}  & = & \big\{z\in\C^k,\ G_i(z)<\infty\big\}\\
          & = &
\big\{z\in\C^k,\mbox{\rm il existe } c>0,\mbox{ \rm tel que }
|f^n(z)|\leq c^{d_i^n}\max\big(|z|^{d_i^n},
1\big)\big\}\\
\K_i & = & \big\{z\in\C^k, \ G_i(z)=0 \big\}
\end{eqnarray*}
et
$\overline\K_i\subset \K_i\cup I_i$.
En particulier, si $f$ est r\'egulier avec $d_m\geq 2$, 
on a $\lambda_1=\lambda_\infty=d_m$, 
$d_t=(d_m)^{l_m-l_{m-1}}\ldots (d_1)^{l_1-l_0}$ et $\K_m=\K$.
\end{theoreme}
\par
On montre le th\'eor\`eme par r\'ecurrence.
Supposons que le th\'eor\`eme et les lemmes suivants soient vrais
jusqu'au rang $i-1$ avec $2\leq i\leq s$. V\'erifions les
au rang $i$. La preuve est aussi valable pour le rang $1$ 
(\voir \'egalement \cite{GuedjSibony}).
\begin{lemme} 
1.  Si $z\in\K_{i-1}$ et $z\rightarrow X_i$, on a  $|f_{(i)}(z)| \sim
|z|^{d_i}$, $|f_{(>i)}(z)|=\o(|z|)^{d_i}$ et $|f(z)|\sim |z|^{d_i}$. 
\par
2. $U_i$ est un ouvert de $\K_{i-1}$, $\K_i$ est un ferm\'e de
$\K_{i-1}$ et $\overline\K_i \subset \K_i\cup I_i$. 
\end{lemme}
\begin{preuve}
1. Puisque $f$ est $i$-r\'egulier, 
$P_{(i)}^+$ ne s'annulle pas sur $X_i$. Par cons\'equent, quand
$z\rightarrow X_i$, on a $|f_{(i)}(z)|\sim |z|^{d_i}$. Puisque $\deg
f_{(>i)}<d_i$, on a $|f_{(>i)}(z)|=\o(|z|^{d_i})$ quand $z\rightarrow
X_i$. 
\par
On sait par hypoth\`ese de r\'ecurrence 
que $\overline\K_{i-1} \subset \K_{i-1}\cup I_{i-1}$  
et que $\K_{i-1}$ est invariant par $f$. 
Par d\'efinition,
$X_i=I_{i-1}\cap \{z_{(>i)}=0\}$. On en d\'eduit que si $z\in \K_{i-1}$
et $z\rightarrow X_i$, on a $f(z)\rightarrow X_i$. D'autre part, 
$f$ \'etant $i$-r\'egulier, on a 
$$I_{i-1}\cap \big\{z_{(\geq i)}=0\big\}=I_{i-1}\cap X_{i-1}=\emptyset.$$ 
Ceci implique en utilisant l'hypoth\`ese de r\'ecurrence 
que lorsque $f(z)\rightarrow I_{i-1}$, on a $|f(z)|\sim
|f_{(\geq i)}(z)|$. Par cons\'equent, si $z\in \K_{i-1}$ et
$z\rightarrow X_i$ on a $|f(z)|\sim |z|^{d_i}$. Observons que
c'est $f_{(\geq i)}$ qui domine sur $\K_{i-1}$ m\^eme si $f_{(1)}$ par
exemple a des termes de degr\'e plus \'elev\'e.
\par
2. D'apr\`es la partie 1, puisque $d_i\geq 2$, 
on peut trouver un voisinage assez petit
$V$ de $X_i$ tel que $f(\K_{i-1}\cap V)\subset V$. Comme $\K_{i-1}$
est invariant, 
pour tout $z\in V$, on a $f^n(z)\rightarrow X_i$. Par cons\'equent,
$\K_{i-1}\cap V\subset U_i$. Par d\'efinition de $U_i$, on a
$U_i=\bigcup_{n\geq 0} f^{-n}(\K_{i-1}\cap V)$. Ceci implique que
$U_i$ est un ouvert de $\K_{i-1}$ et donc $\K_i=\K_{i-1}\setminus U_i$
est un ferm\'e de $\K_{i-1}$. Il reste \`a montrer que $\overline
\K_i\subset \K_i\cup I_i$.
\par
Soit $a\in I_{i-1}\setminus I_i$ et soit
$z\in \K_{i-1}$ tendant vers $a$. On a $f(z)\in \K_{i-1}$ car
$\K_{i-1}$ est invariant. Puisque $a\not\in I_i$, on a
$|f_{(i)}(z)|\sim |z|^{d_i}$. D'autre part,
$|f_{(>i)}(z)|=\o(|z|^{d_i})$ car $\deg f_{(>i)}<d_i$. Par
cons\'equent, $f(z)$ tend vers $X_i$ quand $z\rightarrow a$. 
On conclut que si $z$ est assez
proche de $a$, $f(z)$ appartient \`a $U_i$ et donc $z$ appartient \`a
$U_i$. Ceci implique que $\overline \K_i \subset
\K_i \cup I_i$.
\end{preuve}
{\it Fin de la d\'emonstration du th\'eor\`eme 2.1}--
Montrons d'abord l'existence de la suite $(c_i)$.  
D'apr\`es le lemme 2.2, il existe une
constante $c>0$ telle que pour tout $z\in\K_{i-1}$ on ait
$$|f(z)|\leq c\max\left(|z|^{d_i},1\right).$$
Donc
$$|f^n(z)|\leq c\max\left(|f^{n-1}(z)|^{d_i},1\right)$$
et
$$G_{i,n}(z)\leq \frac{\log c}{d_i^n} +G_{i,n-1}(z).$$
Posons 
$$c_{i,n}:=-\sum_{m\geq n+1} \frac{\log c}{d_i^m}.$$
Il est clair que $c_{i,n}\rightarrow 0$ et que la suite
$G_{i,n}(z)+c_{i,n}$ est d\'ecroissante. Comme les $G_{i,n}(z)$
sont positives,
la limite $G_i(z)$ existe et est positive. 
\par
Montrons les assertions sur les ensembles $\K_{i-1}$ et
$\K_i$.  
Soit $z\not\in
\K_{i-1}$. Il existe $1\leq j\leq i-1$ tel que $z\in U_j$.
D'apr\`es le lemme 2.2, il existe un $c>0$ telle que pour $n$
assez grand  $|f^n(z)|\geq c|z|^{d_j^n}$.
Comme $d_j>d_i$, on a $G_i(z)=+\infty$.
On obtient donc
$$\K_{i-1}   =  \big\{z\in\C^k,\ G_i(z)<\infty\big\}.$$
\par
Soit $z\in \K_i$. Puisque $\K_i$ est invariant, $f(z)\in
\K_i$. D'apr\`es le lemme 2.2, on a
$\overline \K_i\subset \K_i\cup I_i$. Le fait que $I_i\cap
\{z_{(>i)}=0\}=\emptyset$ implique qu'il existe $c'\geq 1$
ind\'ependant de $z$ telle que
$$|f(z)|\leq c'\max\big(|f_{(>i)}(z)|,1\big).$$
Comme $\deg f_{(>i)}=d_{i+1}$, il existe une constante $c''>0$ telle
que
$$|f_{(>i)}(z)|\leq c''\max\left(|z|^{d_{i+1}},1\right).$$
On en d\'eduit que pour une certaine constante $c>0$ on a
$$|f(z)|\leq c\max\left(|z|^{d_{i+1}},1\right).$$
Le fait que $d_{i+1}<d_i$ implique que $G_i(z)=0$. 
On a $\K_i\subset\{z\in\C^k,\ G_i(z)=0\}$.
D'apr\`es le lemme 2.2,
$G_i$ est strictement positive au voisinage de
$X_i$. La relation $G_i\circ f=d_iG_i$ implique que $G_i$ 
est strictement positive sur
$U_i$. On a donc $\K_i\supset\{z\in\C^k,\ G_i(z)=0\}$.
\par
Montrons que $G_i$ est continue et strictement positive sur $U_i$.
D'apr\`es le lemme 2.2, il suffit de le montrer pour $z\in
V\cap\K_{i-1}$ o\`u $V\subset \P^k$ est un voisinage assez petit de $X_i$. 
La positivit\'e est claire. La continuit\'e
r\'esulte de l'in\'egalit\'e:
$$|G_{i,n}(z)-G_{i,n-1}(z)|\leq \frac{\max(\log c,-\log c')}{d_i^n}$$
o\`u $0<c'<c$ sont des constantes telles que $c'|z|^{d_i}\leq
|f(z)|\leq c|z|^{d_i}$ sur $V\cap \K_i$. 
\par
V\'erifions la continuit\'e de $G_i(z)-\log^+|z|$ dans
$U_i\cup I_{i-1}\setminus I_i$. Il suffit de le prouver pour
$z\in V$. Ceci est une cons\'equence de
l'in\'egalit\'e pr\'ec\'edente et de l'\'egalit\'e
$G_{i,0}(z)=\log^+|z|$.
\par
Comme $G_i$ est limite d\'ecroissante d'une suite de fonctions
continues, elle est semi-continue sup\'erieurement sur $\K_{i-1}$. Le fait
qu'elle soit continue, positive sur $U_i$ et nulle sur $\K_i$
implique qu'elle est continue sur $\K_{i-1}=U_i\cup \K_i$.
\par
Dans la suite, on suppose que $f$ est r\'egulier et $d_m\geq 2$.  
On a $X_m=I_{m-1}\not=\emptyset$. Donc
$I_m=\emptyset$. Ceci entra\^{\i}ne que
$$\big\{P_{(1)}^+=\cdots=P_{(m)}^+=0 \big\}=\{z=0\}.$$
Rappelons que $\lambda_n$ et $\lambda_\infty$ sont 
les exposants de Lojasiewicz de $f^n$ et 
l'exposant de Lojasiewicz asymptotique de $f$. 
D'apr\`es le lemme 2.2, on a $\lambda_n\leq d_m^n$. 
Montrons que $\lambda_1\geq d_m$ et $d_t=d$ avec
$d:=(d_1)^{l_1-l_0}\ldots (d_m)^{l_m-l_{m-1}}$. Posons
$\Pi:\C^k\longrightarrow \C^k$ avec
$\Pi(z):=\big(z_{(1)}^{d/d_1},\ldots,z_{(m)}^{d/d_m}\big)$.
C'est une application propre
de degr\'e alg\'ebrique $d/d_m$ et de degr\'e 
topologique $d^{k-1}$. L'application $\Pi\circ f$ est de
degr\'e alg\'ebrique $d$ et se prolonge
holomorphiquement \`a l'infini
car
$$\big\{P_{(1)}^+=\cdots=P_{(m)}^+=0 \big\}=\{z=0\}.$$
On en d\'eduit qu'elle est propre. Son exposant de Lojasiewicz est
\'egal \`a $d$ et
son degr\'e topologique est \'egal \`a $d^k$. Par suite, l'exposant de
Lojasiewicz de $f$ est minor\'e par $d/(d/d_m)=d_m$ et le degr\'e
topologique de $f$ est \'egal \`a $d^k/d^{k-1}=d$. On d\'eduit aussi
que $\lambda_n\geq d_m^n$ et donc $\lambda_1=\lambda_\infty=d_m$.
\par
Puisque $I_m=\emptyset$, l'ensemble $\K_m$ est compact. Il est donc 
\'egal \`a l'ensemble des points d'orbite
born\'ee $\K$.
\par
\hfill $\square$\\
\begin{theoreme} Soit $f:\C^k\longrightarrow\C^k$ comme au
th\'eor\`eme 2.1. Alors pour tout $1\leq r \leq s$ et
$l_{r-1}\leq j <l_r$, on peut d\'efinir le courant $T_j$ de
bidegr\'e $(j,j)$ de $\P^k$ par
$$T_j:=(\ddc G_r)^{j-l_{r-1}+1}\wedge (\ddc
G_{r-1})^{l_{r-1}-l_{r-2}} \wedge \ldots \wedge (\ddc
G_1)^{l_1-l_0}.$$
C'est un courant positif, ferm\'e, de masse $1$, port\'e par
$\overline\K_{r-1}$. Il ne charge pas les ensembles pluripolaires et
on a 
$$f^* T_j=(d_r)^{j-l_{r-1}+1}(d_{r-1})^{l_{r-1}-l_{r-2}}
\ldots (d_1)^{l_1-l_0}T_j.$$
De plus, le courant $T_{l_r-1}$ est port\'e par $\overline\K_r$. 
\end{theoreme}
\begin{preuve} 
D'apr\`es le th\'eor\`eme 2.1, le lemme 2.2 et l'appendice,
le courant $T_j$
est bien d\'efini,
positif, ferm\'e dans $\P^k\setminus I_r$.
Comme $f$ est
$s$-r\'egulier, on a
$$\dim I_r=k-l_r<k-j.$$
\par
Montrons par r\'ecurrence sur $j$ que $T_j$ est de masse 1. Supposons
que c'est
le cas pour $T_{j-1}$.
Posons $\varphi_M(z):=\sup(G_r(z),\log^+|z|-M)$ et $S_M:=\ddc\varphi_M\wedge
T_{j-1}$. La suite $\varphi_M$ d\'ecroit vers $G_r$ sur le support de
$T_{j-1}$. D'apr\`es le lemme 2.2 et la proposition 5.2, $\lim S_M=T_j$
dans $\P^k\setminus I_r$. D'apr\`es la proposition 5.4, la masse de
$S_M$ dans $\P^k$ est \'egale \`a 1. Soit $S$ une valeur adh\'erente
de $(S_M)$ dans $\P^k$. Alors $T_j$ est \'egal \`a $S$ dans
$\P^k\setminus I_r$. Comme $\dim I_r=k-l_r<k-j$, 
les courants $S$ et $T_j$ ne
chargent pas $I_r$. On en d\'eduit que $T_j=S$ et donc 
$T_j$ est de masse 1. 
Le th\'eor\`eme de Skoda \cite{Skoda} entra\^{\i}ne que
$T_j$, qui est de bidegr\'e $(j,j)$, se prolonge en courant invariant,
positif et ferm\'e dans
$\P^k$ qui ne charge pas $I_r$.
D'apr\`es la proposition 5.3, il ne charge pas les
ensembles pluripolaires car la fonction $G_r$ est $T_{j-1}$-p.s.h
pour tout $r\geq 2$ et tout $l_{r-1}\leq j <l_r$.
\par
Pour montrer que $T_j$ est port\'e par $\overline\K_{r-1}$, on peut
supposer que $r\geq 2$. Il suffit de v\'erifier que
$T_{l_{r-1}-1}$ est port\'e par $\overline\K_{r-1}$. On le fait par
r\'ecurrence. Supposons que c'est vrai au rang $r-1$. 
On a
$$T_{l_r-1}=(\ddc G_r)^{l_r-l_{r-1}}\wedge T_{l_{r-1}-1}.$$
Dans $U_r$, on a $G_r=\lim d_{l_{r-1}}^{-n}\log|f^n_{(r)}|$.
D'apr\`es la proposition 5.2, puisque $\big(\ddc
\log|f^{(n)}_{(r)}| \big)^{l_r-l_{r-1}}\wedge T_{l_{r-1}-1}=0$, on a
$(\ddc G_r)^{l_r-l_{r-1}}\wedge T_{l_{r-1}-1}=0$
dans $U_r$. On
en d\'eduit que $T_{l_r-1}$ est port\'e par $\overline\K_r$.
\par
L'application $f$ \'etant ouverte et \`a fibres finies, $f^*$ op\`ere
contin\^ument sur les courants \cite{Meo} et commute avec $\ddc$, ce
qui permet d'obtenir l'\'equation v\'erifi\'ee par $f^*T_j$.
\end{preuve}
\section{Endomorphismes semi-r\'eguliers}
On se propose d'\'etudier d'une
mani\`ere analogue la famille 
des endomorphismes {\it semi-r\'eguliers
ou $(\pi,s)$-r\'eguliers}.
Soient
$1\leq p_1\leq \cdots \leq p_m$ des entiers naturels. Soit
$\pi:\C^k\longrightarrow \C^k$ l'application d\'efinie par
$\pi(z):=\big(\pi_{(1)},\ldots, \pi_{(m)} \big)$ o\`u
$\pi_{(i)}:\C^{l_i-l_{i-1}}\longrightarrow \C^{l_i-l_{i-1}}$ est une
application polynomiale de degr\'e $p_i$ en $z_{(i)}$
qui se prolonge en une application holomorphe de $\P^{l_i-l_j}$ dans
$\P^{l_i-l_j}$. 
Il est clair que $|\pi_{(i)}(z_{(i)})|\sim |z_{(i)}|^{p_i}$ quand
$z_{(i)}$ tend vers l'infini.  On dira qu'une telle application $\pi$
est {\it scind\'ee}.
\par
Soit
$f=(P_1,\ldots,P_k)= \big(P_{(1)},\ldots,P_{(m)}\big)$ un
endomorphisme polynomial ouvert v\'erifiant $\deg P_{(i)} \geq\deg
P_{(i+1)}$. 
Contrairement \`a l'application $\pi$, les degr\'es des fonctions
coordonn\'ees de $f$ sont dans l'ordre d\'ecroissant.
On dit que $f$
est {\it $(\pi,s)$-r\'egulier} 
(resp. {\it $\pi$-r\'egulier}) 
si
l'endomorphisme  $f\circ\pi$ est $s$-r\'egulier (resp.
r\'egulier) o\`u les nombres $1=l_0< l_1< \cdots < l_m=k+1$ utilis\'es 
dans la d\'efinition des applications r\'eguli\`eres sont les m\^emes
que ci-dessus. 
\par
Observons que si $f$ est $(\pi,1)$-r\'egulier, il
est $1$-r\'egulier donc alg\'ebriquement stable. Si on pose
$\pi^+:=\big(\pi^+_{(1)}, \ldots, \pi^+_{(m)}\big)$ o\`u $\pi^+_{(i)}$ est la
partie homog\`ene de plus haut degr\'e de $\pi_{(i)}$, alors $f$ est
$(\pi,s)$-r\'egulier si et seulement s'il est $(\pi^+,s)$-r\'egulier.
Soit $\pi^1$ une application homog\`ene scind\'ee. 
Si $f$ est $(\pi\circ\pi^1,s)$-r\'egulier, alors il est
$(\pi,s)$-r\'egulier.
\par
Il est facile de v\'erifier si un endormorphisme est r\'egulier. Dans
la suite, nous donnons, pour le cas de dimension 2, 
un crit\`ere simple pour savoir si un endomorphisme
polynomial est semi-r\'egulier. L'id\'ee utilis\'ee dans la suite
montre aussi que dans le cas de dimension sup\'erieure \`a 2, on
``devine'' facilement l'application $\pi$ lorsque $f$ est
semi-r\'egulier. 
\par
Consid\'erons dans $\C^2$ l'endomorphisme 
$f(z_1,z_2):=(P_1,P_2)$ avec $d_1:=\deg P_1>d_2:=\deg P_2\geq 1$. 
Cet endomorphisme se prolonge en application m\'eromorphe de $\P^2$
dans $\P^2$. L'ensemble d'ind\'etermination de $f$ est
d\'efini par $I_1:=\{P_1^+=t=0\}$ o\`u $P_i^+$ est la partie
homog\`ene de plus haut degr\'e de $P_i$. On a 
$$X_1:=\{t=0\}\cap\{z_{(>1)}=0\}=f(\{t=0\}\setminus I_1)=[1:0:0].$$
L'endomorphisme $f$ est
$(\pi,1)$-r\'egulier pour une application scind\'ee $\pi$, 
s'il est $1$-r\'egulier. Dans le 
cadre consid\'er\'e, cela \'equivaut \`a dire qu'il
est alg\'ebriquement
stable, \cad $X_1\cap I_1=\emptyset$. Autrement dit, 
le coefficient de $z_1^{d_1}$ dans $P_1$ est non nul. 
\par
Dans la suite, on suppose que $f$ est alg\'ebriquement stable. 
Notons $\Sigma_i\subset \N^2\subset \R^2$ 
l'ensemble des couples $(m,n)$ tels que le
coefficient de $z_1^mz_2^n$ dans $P_i$ soit non nul. 
Puisque $\deg P_2< \deg P_1 =d_1$,  
les $\Sigma_i$ se trouvent au dessous de la droite
$m+n=d_1$. L'endomorphisme $f$ \'etant alg\'ebriquement stable, $(d_1,0)\in
\Sigma_1$. Soit $\D$ une droite dans $\N^2$. On note $P_i^{\D}$ la
somme des termes $z_1^mz_2^n$ dans $P_i$ avec $(m,n)\in \D$. Si $\D$
est d\'efinie par l'\'equation $pm+qn+r=0$, on dira que $-p/q$ est 
{\it la pente}
de cette droite. 
\par
On note $\D_1$ la droite v\'erifiant les propri\'et\'es suivantes:
\begin{enumerate}
\item $\D_1$ passe par $(d_1,0)$ et au moins un autre point de
  $\Sigma_1\cup\Sigma_2$. 
\item L'ensemble $\Sigma_1\cup \Sigma_2$ est au dessous de $\D_1$.
\end{enumerate}
La droite $\D_1$ est celle de pente maximale, passant par $(d_1,0)$ et
v\'erifiant la condition 2. Cette pente est comprise entre $-1$ et $0$.
Notons $\D_2$ la droite parall\`ele \`a $\D_1$ et
v\'erifiant:
\begin{enumerate}
\item $\D_2$ passe par au moins un point de
  $\Sigma_2$. 
\item L'ensemble $\Sigma_2$ est au dessous de $\D_2$.
\end{enumerate}
Il est clair que $\D_2$ est au dessous de $\D_1$.
On a la proposition suivante:
\begin{proposition} Soit $f$ un endomorphisme de $\C^2$
  alg\'ebriquement stable comme ci-dessus. Alors $f$
  est semi-r\'egulier si et seulement si la pente de $\D_1$ est non
  nulle et l'ensemble
$\big\{z\in\C^2,\ P_1^{\D_1}(z)=P_2^{\D_2}(z)=0 \big\}$ est r\'eduit \`a
  $\{0\}$.
Dans ce cas, 
si la pente de $\D_1$ est \'egale \`a $-p_1/p_2$ avec $p_1\in \N^*$ et
$p_2\in \N^*$, $f$ est $\pi$-r\'egulier pour
$\pi(z):=(z_1^{p_1},z_2^{p_2})$. 
\end{proposition}
\begin{preuve}
Supposons que $\big\{z\in\C^2,\ P_1^{\D_1}(z)
=P_2^{\D_2}(z)=0 \big\}=\{0\}$. Puisque
$\Sigma_1\cup\Sigma_2$ est au dessous de la droite de pente $-1$
passant par $(d_1,0)$, la pente de $\D_1$ est plus grande ou \'egale
\`a $-1$. Soit $-p_1/p_2$ la pente de $\D_1$ o\`u  
$p_1$ et $p_2$ sont des entiers positifs. On a $p_1\leq p_2$. Posons
$\pi(z):=(z_1^{p_1},z_2^{p_2})$. On a $P^{\pi+}_i=P_i^{\D_i}\circ
\pi$. On en d\'eduit que $\big\{z\in\C^2,\
P^{\pi+}_1(z)=P^{\pi+}_2(z)=0 \big\}=\{0\}$. Comme $\D_2$ est au dessous de
$\D_1$, on a $\deg P^{\pi+}_1\geq \deg P^{\pi+}_2$. On a montr\'e que
$f\circ\pi$ est r\'egulier. Donc $f$ est $\pi$-r\'egulier.
\par
Supposons maintenant que $f$ est $\pi$-r\'egulier avec
$\pi(z)=(z_1^{p_1},z_2^{p_2})$ et $p_1\leq p_2$. 
On note $\D_i'$ la droite de pente
$-p_1/p_2$ telle que
\begin{enumerate}
\item $\D_i'$ passe par au moins un point de $\Sigma_i$.
\item $\Sigma_i$ est au dessous de $\D_i'$. 
\end{enumerate}
On a $P^{\pi+}_i=P^{\D_i'}_i\circ \pi$ et donc 
$\big\{z\in\C^2,\ P^{\D_1'}_1(z)=P^{\D_2'}_2(z)=0 \big\}=\{0\}$.
Puisque $\deg P^{\pi+}_1\geq
\deg P^{\pi+}_2$, la droite $\D_2'$ et l'ensemble $\Sigma_1\cup
\Sigma_2$ sont au dessous de $\D_1'$.
\par
Comme $f\circ\pi$ est
r\'egulier, $P^{\pi+}_1$ ne s'annulle pas en $X_1$. Par cons\'equent,
$P^{\D_1'}_1$ contient un mon\^ome en $z_1$ et donc $\D_1'$ passe
par $(d_1,0)$. On en d\'eduit que la pente de $\D_1'$ est plus petite
que celle de $\D_1$.
Si la pente de $\D_1'$ est \'egale \`a celle de $\D_1$, alors
$\D_i'=\D_i$, donc  
$$\big\{z\in\C^2,\
P^{\D_i}_1(z)=P^{\D_i}_2(z)=0 \big\}=\{0\}.$$
\par
Il reste \`a consid\'erer le cas o\`u 
la pente de $\D_1'$ est strictement plus petite que celle de
$\D_1$. Dans ce cas, $\D_1'$ ne passe par aucun point de
$\Sigma_1$ sauf le point $(d_1,0)$. Par cons\'equent, $P_1^{\D_1'}$
est un mon\^ome en $z_1$. Ceci implique que  $P_2^{\D_2'}$ contient
un mon\^ome en $z_2$ et donc $\D_2'$ passe par un point $(0,d)$.  
Comme la pente $\D_2$ est strictement plus
grande que celle de $\D_2'$, par d\'efinition, $\D_2$ ne passe par
aucun point de $\Sigma_2$ sauf le point $(0,d)$. On conclut que
$P^{\D_2}_2$ est un mon\^ome en $z_2$. Le fait que $P^{\D_1}_1$ contient
un mon\^ome en $z_1$ implique que 
$$\big\{z\in\C^2,\
P^{\D_1}_1(z)=P^{\D_2}_2(z)=0 \big\}=\{0\}.$$
\end{preuve}
\par
La proposition 3.1 permet de v\'erifier facilement dans le cas de
dimension 2 si un endomorphisme est semi-r\'egulier. Donnons des
exemples. 
\begin{exemple} \rm
Consid\'erons 
les endomorphismes 
$$f(z):=(z_1^6-z_2^4,z_1^3-2z_2^2+z_2) \ \ 
\mbox{ et } \ \ 
g(z):=(z_1^6-z_2^4,z_1^3-z_2^2+z_2).$$ 
Dans les deux cas, on a
$\D_1=\{2m+3n=12\}$ et $\D_2=\{2m+3n=6\}$. Dans le premier cas,
on a $P_1^{\D_1}(z)=z_1^6-z_2^4$ et  $P_2^{\D_1}(z)=z_1^3-2z_2^2$;
la condition de la proposition 3.1 est v\'erifi\'ee.
Dans le second
cas, on a $P_1^{\D_1}(z)=z_1^6-z_2^4$ et
$P_2^{\D_1}(z)=z_1^3-z_2^2$;
la condition de la proposition 3.1  
n'est pas v\'erifi\'ee. L'endomorphisme $f$ est $\pi$-r\'egulier
pour $\pi(z):=(z_1^2,z_2^3)$. L'endomorphisme $g$ n'est pas
semi-r\'egulier. Puisque $f$ est 1-r\'egulier, on peut d\'efinir 
$X_1=[1:0:0]$, $I_1=[0:1:0]$, $U_1$ le bassin de $X_1$ et
$\K_1=\C^2\setminus U_1$
v\'erifiant $\overline \K_1=\K_1\cup I_1$. Si $z\in \K_1$ tend vers
$I_1$, $f(z)$ tend aussi vers $I_1$. Par cons\'equent, sa deuxi\`eme
coordonn\'ee domine la premi\`ere; on a
$z_1^6-z_2^4=\o(z_1^3-2z_2^2+z_2)$. On en d\'eduit que
$z_1^6\simeq z_2^4$ et donc $z_1^3\simeq z_2^2$ ou $z_1^3\simeq
-z_2^2$. Cela entra\^{\i}ne que la deuxi\`eme composante de $f(z)$
v\'erifie $|z_1^3-2z_2^2+z_2|\sim |z_2|^2\sim |z|^2$. Cette derni\`ere
relation permet de
d\'efinire la deuxi\`eme fonction de Green comme dans le cas des
applications r\'eguli\`eres.
\end{exemple}
\par 
Etudions maintenant la dynamique d'une application
$(\pi,s)$-r\'eguli\`ere $f$ g\'en\'erale. 
On pose $P^\pi_{(i)}:=f_{(i)}\circ\pi$ et on note $P^{\pi+}_{(i)}$
sa partie homog\`ene de plus haut degr\'e $d^\pi_i$.
Soient $\check{P}_{(i)}$ les polyn\^omes constitu\'es par certains 
mon\^omes de $P_{(i)}$ avec les m\^emes coefficients et v\'erifiant la
relation
$P^{\pi +}_{(i)}=\check{P}_{(i)}\circ \pi^+$. 
Cette derni\`ere relation n'assure pas l'unicit\'e de
$\check{P}_{(i)}$. On garde les m\^emes coefficients pour l'assurer.
Notons 
$\check{P}^+_{(i)}$ la partie homog\`ene  de degr\'e
$\alpha_i:=d_i^\pi/p_i$
de  $\check{P}_{(i)}$. En particulier, si $\alpha_i$ n'est pas entier,
on a $\check{P}^+_{(i)}=0$. 
\par
Posons $I_0=\{t=0\}$, $X_0=\emptyset$. 
Nous allons d\'efinir d'une mani\`ere analogue au cas
$s$-r\'egulier, les autres ensembles $X_i$ et $I_i$. 
Lorsque $I_{i-1}$ est d\'efini, on pose toujours $X_i:=I_{i-1}\cap
\{z_{(>i)}=0\}$. En particulier, on a $X_1=\{t=0\}\cap \{z_{(>1)}=0\}$. 
Il reste \`a d\'efinir les ensembles $I_i$.
Les notations
\'etant assez compliqu\'ees, nous expliquons d'abord des cas
simples. Si $p_1=\cdots =p_m$, l'application $f$ est
$s$-r\'eguli\`ere; les ensembles $X_i$ et $I_i$ sont d\'efinis
exactement de la m\^eme mani\`ere qu'au paragraphe
pr\'ec\'edent. Si $p_1<\cdots<p_m$ (\voir l'exemple 3.2), 
le fait que $f\circ\pi$ soit
$s$-r\'egulier implique que
$P^+_{(1)}=\check{P}_{(1)}^+$ 
ne d\'epend pas de
$z_{(>1)}$. Il est clair qu'il faut prendre 
$$I_1:=I_0\cap\{z_{(1)}=0\}=I_0\cap\{\check{P}_{(1)}^+=0\}.$$
La croissance des $p_i$ implique aussi
que les termes de $\check{P}^{+}_{(2)}$, qui sont
ind\'ependants de $z_{(1)}$, ne d\'ependent que de $z_{(2)}$. 
Rappelons que $I_2$ appara\^{\i}t comme l'ensemble d'ind\'etermination
d'une restriction convenable de $f$ \`a $I_1=I_0\cap\{z_{(1)}=0\}$. 
Il est
donc aussi clair qu'il faut poser  
$$I_2:=I_1\cap\{z_{(2)}=0\}=I_1\cap \{\check{P}^+_{(2)}=0\}.$$ 
Le m\^eme raisonnement est valable pour
$I_i$ avec $1\leq i\leq s$. 
\par
Pour le cas g\'en\'eral, la
d\'efinition des $X_i$ et $I_i$ tient compte des 
deux cas particuliers ci-dessus.
Posons pour tout $1\leq i\leq m$
$$X_i:=I_{i-1}\cap \big\{z_{(> i)}=0 \big\}$$
et
$$I_i:=I_{i-1}\cap
\big\{\check{P}^+_{(i)}=0 \big\}\cap 
\big\{z_{(j)}=0 \mbox{ pour tout } j \mbox{ tel
  que } p_j<p_i \big\}.$$
\begin{lemme}
Si $f$ est $(\pi,s)$-r\'egulier, alors $X_i\cap I_i=\emptyset$
pour tout $1\leq i\leq s$. De plus, on a $\dim X_i=l_i-l_{i-1}-1$ et
$\dim I_i =k-l_i$. 
\end{lemme}
\begin{preuve} Soit $0\leq j <i$ l'entier minimal tel que
  $p_{j+1}=p_i$. Posons $I^1_0=\P^k$, 
$$I_i^1:=\big\{ z_{(\leq j)} =0 \big\} \cap \big\{ \check{P}^+_{(\leq
  i)} =0 \big\}$$
et
$$X^1_i:= I^1_{i-1}\cap \big\{ z_{(>i)}=0 \big\}.$$
On a $X_i=X_i^1\cap\{t=0\}$ et $I_i=I_i^1\cap\{t=0\}$. 
Il faut montrer que $X_i^1\cap I_i^1\cap \C^k=\{0\}$ et que 
$\dim X_i^1=l_i-l_{i-1}$, $\dim I_i^1= k-l_i+1$. Posons
$X_i^2:=(\pi^+)^{-1}(X_i^1\cap \C^k)$ et 
$I_i^2:=(\pi^+)^{-1}(I_i^1\cap \C^k)$. 
En utilisant la d\'ecroissance des degr\'es de $f$
et la croissance des degr\'es de $\pi$, on v\'erifie
qu'aucun terme de $P^{\pi +}_{(i)} - \check{P}_{(i)}^+ \circ
\pi^+$ n'est ind\'ependant de $z_{(\leq j)}$.
On a
\begin{eqnarray*}
X_i^2\cap I_i^2 & = & \big\{ z_{(\leq j)}=0 \big\} \cap
\big\{\check{P}^+_{(\leq i)} \circ \pi^+ =0 \big\} \cap \{z_{(>i)}=0\}
\cap \C^k
  \\
& \subset & \big\{ P^{\pi +}_{(\leq i)} =0 \big\} \cap \big\{ z_{(>i)}=0
\big\} \cap \C^k =\{0\}.
\end{eqnarray*}
La derni\`ere intersection est r\'eduite \`a $\{0\}$ car $f\circ\pi^+$
est $s$-r\'egulier. On d\'eduit de la propri\'et\'e pr\'ec\'edente que
$X_i\cap I_i=\emptyset$.
\par
On a 
\begin{eqnarray*}
I^2_i & = & \big\{ z_{(\leq j)} =0 \big\} \cap \big\{
\check{P}^+_{(\leq i)} \circ \pi^+ =0\big\} \cap \C^k \\
& = & \big\{ z_{(\leq j)} =0 \big\} \cap \big\{
\check{P}^+_{(r)} \circ \pi^+ =0 \mbox{ pour tout } j+1\leq r\leq i 
\big\} \cap \C^k.
\end{eqnarray*}
Par cons\'equent, on a $\dim I_i^2\geq k-l_i+1$. D'autre part,
$X_i^2=I_{i-1}^2\cap \big\{ z_{(>i)} =0 \big\}$. Donc $\dim X_i^2 \geq
l_i-l_{i-1}$. Le fait que $X_i^2\cap I_i^2=\{0\}$ implique que $\dim
X_i^2 = l_i-l_{i-1}$ et $\dim I_i^2 = k-l_i+1$. On en d\'eduit que 
$\dim X_i=l_i-l_{i-1}-1$ et $\dim I_i=k-l_i$. 
\end{preuve}
\par
Posons $G_0:=0$, $\K_0:=\C^k$. 
Lorsque les expressions suivantes ont un sens, on pose
$$G_{i,n}(z):=\frac{\log^+|f^n(z)|}{\alpha_{i}^n}$$
$$G_i(z):=\lim_{n\rightarrow\infty} G_{i,n}(z)$$
$$U_i:=\big\{z\in\C^k,\  f^n(z) \mbox{ tend vers } X_i \big\}$$
et
$$\K_i:=\K_{i-1}\setminus U_i$$
pour tout $1\leq i\leq m$.
L'ensemble $U_i$ est le {\it
bassin d'attraction}
de $X_i$. Il est clair que $G_i\circ f=
\alpha_{i}G_i$, $f^{-1}(U_i)\subset f(U_i)=U_i$ et
$f^{-1}(\K_i)\subset f(\K_i)=\K_i$. Observons ici que la suite 
des $\alpha_i$ est d\'ecroissante.
On obtient
des r\'esultats analogues au cas r\'egulier.
\begin{theoreme} Soit $f:\C^k\longrightarrow \C^k$ un
endomorphisme polynomial $(\pi,s)$-r\'egulier comme
pr\'ec\'edemment. On suppose que
$1\leq s\leq m-1$ ou bien que $s=m$ et $\alpha_m>1$.
Alors
pour tout $1\leq i\leq s$, il existe une suite de nombres r\'eels
positifs 
$(c_{i,n})\rightarrow 0$ telle que la suite de fonctions
$G_{i,n}+c_{i,n}$
d\'ecro\^{\i}t sur $\K_{i-1}$ vers une fonction $G_i$ continue,
invariante: $G_i\circ f=\alpha_{i} G_i$. De plus, la fonction 
$G_i(z)-\log^+|z|$ est
continue sur $\overline\K_{i-1}\setminus I_i$ et on a 
\begin{eqnarray*}
\K_{i-1}  & = & \big\{z\in\C^k,\ G_i(z)<\infty \big\}\\
          & = &
\big\{z\in\C^k,\mbox{\rm\ il existe } c>0 \mbox{ \rm tel que }
|f^n(z)|\leq c^{\alpha_i^n}\max\big(|z|^{\alpha_{i}^n},
1\big)\big\}\\
\K_i & = &\big\{z\in\C^k, \ G_i(z)=0\big\}
\end{eqnarray*}
et
$\overline\K_i\subset \K_i\cup I_i$.
En particulier, si $f$ est $\pi$-r\'egulier et $\alpha_m>1$, 
on a $\lambda_1=\lambda_\infty=\alpha_m$,  
$d_t=(\alpha_m)^{l_m-l_{m-1}}\ldots (\alpha_1)^{l_1-l_0}$ et $\K_m=\K$.
\end{theoreme}
\begin{theoreme} Soit $f:\C^k\longrightarrow\C^k$ comme au
th\'eor\`eme 3.4. Alors pour tout $1\leq i \leq s$ et
$l_{i-1}\leq j <l_i$, on peut d\'efinir le courant $T_j$ de
bidegr\'e $(j,j)$ de $\P^k$ par
$$T_j:=(\ddc G_i)^{j-l_{i-1}+1}\wedge (\ddc
G_{i-1})^{l_{i-1}-l_{i-2}} \wedge \ldots \wedge (\ddc
G_1)^{l_1-l_0}.$$
C'est un courant positif, ferm\'e, de masse $1$, port\'e par
$\overline\K_{i-1}$. Il ne charge pas les ensembles pluripolaires et 
on a
$$f^* T_j=(\alpha_i)^{j-l_{i-1}+1}
(\alpha_{i-1})^{l_{i-1}-l_{i-2}}
\ldots (\alpha_1)^{l_1-l_0}T_j.$$
De plus, le courant $T_{l_i-1}$ est port\'e par $\overline\K_i$. 
\end{theoreme}
\par
Pour d\'emontrer ces th\'eor\`emes, nous allons montrer des
in\'egalit\'es analogues que celles du lemme 2.2. Puisque $f\circ\pi$
est r\'egulier, il est plus facile d'utiliser $\pi(z)$ au lieu de $z$ 
comme ``coordonn\'ees''. L'application $\pi$ n'est pas inversible en
g\'en\'eral mais elle est scind\'ee. Ceci nous permet de travailler
avec $\pi^{-1}$ comme avec une application polynomiale. 
\par
Posons $I_0^\pi:=\{t=0\}$, $X_0^\pi=\emptyset$ 
et pour tout $1\leq i\leq m$
$$X_i^{\pi}:=I_{i-1}^{\pi}\cap \big\{z_{(> i)}=0 \big\}$$
et
$$I_i^{\pi}:=I_{i-1}^{\pi}\cap \big\{P^{\pi+}_{(i)}=0 \big\}.$$
On a $X^\pi_i\cap I^\pi_i=\emptyset$ pour $1\leq i \leq s$.
Soit $\K^\pi_0=\C^k$. Posons
$$U^\pi_i:=\big\{z\in\C^k,\  \pi^{-1}\circ f^n\circ\pi(z) 
\mbox{ tend vers } X^\pi_i\big\}$$
et
$$\K^\pi_i:=\K^\pi_{i-1}\setminus U^\pi_i$$
pour tout $1\leq i\leq s$. 
Observons que $\pi^{-1}\circ f^n\circ\pi(z)$ contient plusieurs
points dont les modules sont comparables quand $f^n\circ\pi(z)$ tend
vers l'infini.
Ceci est en fait une cons\'equence de la propri\'et\'e ``$\pi$ est
scind\'ee''. 
Par d\'efinition, 
$\K^\pi_i$ et $U^\pi_i$ sont invariants par $\pi^{-1}\circ
f\circ\pi$ et par $\pi^{-1}\circ f^{-1}\circ\pi$. 
\begin{lemme} 
1. Pour $z\in \K^\pi_{i-1}$ tendant vers $X^\pi_i$, on a
$$|f_{(i)}\circ\pi(z)|\sim |z|^{d^\pi_i} \mbox{ et } 
|f_{(>i)}\circ\pi(z)| = \o(|z|^{d^\pi_i}).$$
\par
2. Pour tout $z\in\K^\pi_{i-1}$ tendant vers $X^\pi_i$, on a
$|\pi^{-1}\circ f\circ\pi(z)|\sim |z|^{\alpha_i}$ et $|\pi^{-1}\circ
f\circ \pi(z)|_{(>i)}=\o(|z|^{\alpha_i})$.
\par
3. $U^\pi_i$ est un ouvert de $\K^\pi_{i-1}$, $\K^\pi_i$ est un
ferm\'e de $\K^\pi_{i-1}$ et
$\overline\K^\pi_i\subset \K^\pi_i\cup I^\pi_i$.
\par 
4. Pour $z\in\K^\pi_i$ tendant vers $I^\pi_i$, on a
$|P^\pi_{(i)}(z)|=\o(|z|^{d_i^\pi})$.  
\end{lemme}
\begin{preuve}
On montre le lemme par r\'ecurrence. 
On suppose
qu'il est vrai jusqu'au rang $i-1$.
\par
1. La partie homog\`ene de plus haut degr\'e $P^\pi_{(i)}$ de
$f\circ\pi$ ne s'annulle pas sur $X^\pi_i$ car $f\circ\pi$ est
$s$-r\'egulier. Par cons\'equent, pour $z\in \K^\pi_{i-1}$ tendant vers 
$X^\pi_i$, on a
$|f_{(i)}\circ\pi(z)| \sim |z|^{d^\pi_i}$. Quant \`a
$f_{(>i)}\circ\pi$, il est de degr\'e $d^\pi_{i+1}<d^\pi_i$. On a donc 
$$|f_{(>i)}\circ\pi(z)| = \o\big(|z|^{d^\pi_i}\big).$$
\par
2. Par hypoth\`ese de r\'ecurrence, $\overline\K^\pi_{i-1}
\subset \K^\pi_{i-1}\cup I^\pi_{i-1}$. D'autre part, 
$$I^\pi_{i-1}\cap \{z_{(\geq i)}=0\}=\emptyset.$$
On en d\'eduit que pour $z\in\K^\pi_{i-1}$
tendant vers $I^\pi_{i-1}$, on a 
$|z|\sim |z|_{(\geq i)}$. Comme $\pi^{-1}\circ f\circ \pi$
pr\'eserve $\K^\pi_{i-1}$, on a
$$|\pi^{-1}\circ f\circ\pi(z)|\sim |\pi^{-1}\circ f\circ\pi(z)|_{(\geq
  i)}.$$
D'apr\`es la partie 1, si $z\in\K^\pi_{i-1}$ tend vers $X^\pi_i$, 
la derni\`ere relation et la croissance des $p_i$ impliquent  
$$|\pi^{-1}\circ f\circ\pi(z)|\sim |\pi^{-1}\circ f\circ\pi(z)|_{(i)}
\sim |f_{(i)}\circ\pi(z)|^{1/p_i} \sim |z|^{\alpha_i}$$
et
$$|\pi^{-1}\circ f\circ\pi(z)|_{(>i)}=\o(|z|^{d^\pi_i/p_i})=
\o(|z|^{\alpha_i}).$$ 
\par
3. Soit $V$ un voisinage suffisamment petit de $X^\pi_i$. La partie 2
implique $V\cap \K^\pi_{i-1} \subset U^\pi_i$. 
Par
d\'efinition, on a $U^\pi_i:=\bigcup_{n\geq 0} \pi^{-1}\circ f^{-n}
\circ \pi (\K^\pi_{i-1} \cap V)$. C'est donc un ouvert de
$I^\pi_{i-1}$. Par suite, $\K^\pi_i$ est un ferm\'e de $I^\pi_{i-1}$.
\par
Fixons un voisinage $W$ de $I^\pi_i$. 
Observons que lorsque
$z\in\K^\pi_{i-1}$ tend vers
$I^\pi_{i-1}\setminus W$, la partie 1 du lemme est encore vraie; par
suite, $\pi^{-1}\circ f\circ\pi(z)$ tend vers $X^\pi_i$. Par
cons\'equent, pour $z\in\K^\pi_{i-1}$ suffisamment proche de
$I^\pi_{i-1}\setminus W$, on a $z\in U^\pi_i$. 
D'o\`u $\overline\K^\pi_i\subset \K^\pi_i\cup
I^\pi_i$. 
\par
4. Soit $(z^{(n)})\subset \K^\pi_{i-1}$ une suite tendant vers
$I^\pi_{i-1}$. Si $|P^\pi_{(i)}(z^{(n)})|\sim
|z^{(n)}|^{d^\pi_i}$, alors comme dans la partie 1, on montre que
$$|f_{(i)}\circ\pi(z^{(n)})|\sim |z^{(n)}|^{d^\pi_i} \mbox{ et } 
|f_{(>i)}\circ\pi(z^{(n)})| = \o(|z^{(n)}|^{d^\pi_i}).$$
Par cons\'equent, $\pi^{-1}\circ f\circ\pi(z^{(n)})$ tend vers
$X^\pi_i$. D'apr\`es la partie 3, $z^{(n)}$ appartient \`a $U^\pi_i$
pour $n$ assez grand. Ceci d\'emontre la partie 4.
\end{preuve}
\begin{lemme} 
1. Si $w^*\in\K_{i-1}$ tend vers $X_i$, on a
$|f(w^*)|\sim |w^*|^{\alpha_i}$.
\par
2. On a $\pi(\K^\pi_i)=\K_i$, $\pi(U^\pi_i)=U_i$ et 
$\overline\K_i\subset \K_i\cup I_i$.
\end{lemme}
\begin{preuve} On montre le lemme par r\'ecurrence. Supposons qu'il
  est vrai jusqu'au rang $i-1$.
\par
1. Soit $w^*=\pi(z^*)\in\K_{i-1}$ tendant vers $X_i$. Par hypoth\`ese de
r\'ecurrence, on peut choisit 
$z^*\in\K^\pi_{i-1}$ et $z^*$ tend vers $I^\pi_{i-1}$ quand 
$w^*$ tend vers $X_i$. Puisque
les degr\'es $p_i$ des composantes $\pi_{(i)}$ de $\pi$ sont 
croissantes, $z^*$ tend vers $X^\pi_i$.
D'apr\`es le lemme 3.3, on a
\begin{eqnarray*}
X_i\cap \big\{ w_{(i)}=0 \big\} & = & I_{i-1} \cap \big\{w_{(\geq i)}
=0\big\} \\
& = & I_{i-1}\cap X_{i-1} =\emptyset.
\end{eqnarray*}
De plus, $X_i\subset \big\{ w_{(>i)}=0 \big\}$ et la suite $(p_i)$ est
croissante. On en d\'eduit que
$$|z^*|\sim |z^*_{(i)}|\sim |w^*_{(i)}|^{1/p_i} \sim |w^*|^{1/p_i}.$$
La premi\`ere relation de la derni\`ere ligne vient du fait que
$z^*\rightarrow X_i^\pi$.
\par
L'invariance de $\K_{i-1}$ implique que 
$f(w^*)$ appartient \`a $\K_{i-1}$ et 
tend vers $I_{i-1}$. De plus,
$I_{i-1} \cap \{w_{(\geq i)}=0\}=\emptyset$. On en d\'eduit 
$$|f(w^*)|  \sim  |f_{(\geq i)}(w^*)| = |f_{(\geq i)}\circ\pi(z^*)|.$$
D'apr\`es le lemme 3.6, ceci implique
$$|f(w^*)| \sim |f_{(i)}\circ\pi(z^*)| \sim  |z^*|^{d^\pi_i}  
\sim |w^*|^{d^\pi_i/p_i}= |w^*|^{\alpha_i}.$$
\par
2. Soit $w^*=\pi(z^*)\in U_i$. Puisque $f^n(w^*)$ tend vers $X_i$ et
que la suite $(p_i)$ est croissante,
$\pi^{-1}\circ f^n(w^*)=\pi^{-1}\circ f^n\circ \pi(z^*)$ 
tend vers $X^\pi_i$. Par cons\'equent, $z^*\in U^\pi_i$. On a donc
$U_i\subset \pi(U^\pi_i)$.  
\par
Soit maintenant $z^*\in U^\pi_{i-1}$. 
On a $f^n\circ \pi(z^*) =f\circ \pi \circ (\pi^{-1}\circ f\circ\pi)^{n-1}
(z^*)$. 
On applique les parties
1 et 2 du lemme 3.6 aux points de $(\pi^{-1}\circ
f\circ \pi)^{n-1}(z^*)$ qui tendent vers $X_i^\pi$. 
On obtient que la suite
$(f^n\circ\pi(z^*))$ tend vers $\K_{i-1}\cap \{w_{(>i)}=0\}=X_i$. Par
cons\'equent, $\pi(z^*)\in U_i$ et $\pi(U^\pi_i)\subset
U_i$. On a montr\'e $\pi(U^\pi_i)=U_i$. On en d\'eduit que
$\K_i=\pi(\K^\pi_i)$.  
\par
Soit $w^{(n)}=\pi(z^{(n)}) \subset \K_i$ 
une suite tendant vers $a\in I_{i-1}$ avec $z^{(n)}\in \K_i^\pi$. 
Puisque $z^{(n)}$ tend vers $I_i^\pi$, on a $|z^{(n)}|\sim
|z^{(n)}|_{(\geq i)}$. Par cons\'equent, si $p_j<p_i$, on a
$a_{(j)}=0$, \cad que $a\in\{z_{(j)}=0\}$. 
\par
D'apr\`es la partie 4 du lemme 3.6, on a 
$$\check P_{(i)}(w^{(n)}) = \o\big(|z^{(n)}|^{d_i^\pi}\big) =
\o\big(|z^{(n)}_{(\geq i)}|^{d_i^\pi}\big) = 
\o\big(|w^{(n)}_{(\geq i)}|^{\alpha_i}\big).$$
Ceci implique que $\check P_{(i)}^+(a)=0$. 
\par
On a montr\'e que $w^{(n)}$ tend vers
$a\in I_i$. Donc $\overline\K_i\subset
\K_i\cup I_i$.
\end{preuve}
\par
Utilisant ces deux derniers lemmes, 
on montre les th\'eor\`emes 3.4 et 3.5 de m\^eme
mani\`ere que dans le cas des applications r\'eguli\`eres.
\par
\hfill $\square$
\\
\begin{proposition} Soit $f$ un endomorphisme $\pi$-r\'egulier avec
  $\alpha_m>1$, comme
  pr\'ec\'edemment et soit $\K$ l'ensemble des points d'orbite born\'e. 
Posons $M:=\lim \sup_\K
  \|\mbox{\rm D}f^n\|^{1/n}$ 
o\`u ${\rm D}$ d\'esigne la d\'eriv\'ee. 
Alors 
\begin{enumerate}
\item Pour tous $1\leq i\leq m$ et $0<a_i<\log \alpha_i/\log M$,
il existe une constante $c>0$ telle que 
si  $z\in\K_{i-1}$ on ait $G_i(z)\leq c \delta(z)^{a_i}$ o\`u
  $\delta(z)$ d\'esigne la distance entre $z$ et $\K$. 
\item La mesure $\mu:=T_k$ ne charge pas les ensembles de
  dimension de Hausdorff 
$a$ pour tout $a <\log d_t/\log M$. 
\end{enumerate}
\end{proposition}
\begin{preuve} 1. Puisque la fonction $G_i$ est \`a croissance
  logarithmique, il suffit de montrer que
$G_i(z)\leq c \delta(z)^{a_i}$
  dans un voisinage fixe de $\K$. 
Soit $W$ un voisinage assez petit de $\K$ et soit 
  $N$ assez grand tels que $a_i <a_i^*:=\log \alpha_i/\log M_N$ pour
  tout $1\leq i\leq m$ o\`u $M_N:=\sup_W \|\mbox{D} f^N\|^{1/N}$.
Notons
  $\delta>0$ la distance entre $\K$ et $\partial W$. On choisit une
  constante $A>0$ telle que $G_i(z)\leq A$ pour tout $z\in \K_{i-1}\cap W$ et
  tout $1\leq i\leq m$.
\par
Il suffit de montrer
 qu'il existe $c>0$ tel
  que pour tout $z\in \K_{i-1}\cap W$ on a $|G_i(z)|\leq c
  \delta(z)^{a_i^*}$. Comme $G_i=0$ sur $\K$, il suffit de
  consid\'erer le cas $z\not\in \K$. 
Soit $n$ l'entier minimal tel que 
  $f^{Nn}(z)\not\in W$. Si $x\in \K$ est un point tel que
  $|x-z|=\delta(z)$ alors on a
$$\delta\leq |f^{Nn}(x)-f^{Nn}(z)| \leq (M_N)^{Nn} |x-z|= (M_N)^{Nn} 
\delta(z)$$
donc
$$(M_N)^{Nn}\geq \frac{\delta}{\delta(z)}\ .$$
Cela entra\^{\i}ne que pour une certaine constante $c'>0$ on a
$$\frac{1}{\alpha_i^{N(n-1)}}=\alpha_i^NM_N^{-Nna_i^*} 
\leq 
\alpha_i^N \left(\frac{\delta(z)}{\delta}\right)^{a_i^*} 
\leq c'\delta(z)^{a_i^*}.$$   
Puisque $f^{N(n-1)}(z)\in \K_{i-1}\cap W$, en posant $c=Ac'$, on a 
$$G_i(z)=\frac{G_i(f^{N(n-1)}(z))}{\alpha_i^{N(n-1)}}\leq
Ac'\delta(z)^{a_i^*} =c\delta(z)^{a_i^*}.$$ 
\par
2. Comme dans la partie 1, on peut choisir $W$ de sorte que 
$a <a^*:=\log d_t/\log M_N$. Puisque $d_t=\alpha_m^{l_m-l_{m-1}}\ldots
  \alpha_1^{l_1-l_0}$, on a $a^*=(l_m-l_{m-1})a^*_m+\cdots +
  (l_1-l_0)a^*_1$.
Soit $\Sigma\subset \K$ un ensemble de dimension de Hausdorff $a$. 
Fixons un $\epsilon>0$ assez petit. Comme $a<a^*$, pour 
$r>0$ assez petit on
peut recouvrir $\Sigma$ par $\epsilon r^{-a^*}$  boules $B_n(r)$ de 
centres $x_n\in \K$ et de rayon $r$. Soit $\chi_n$ une fonction
positive \`a support dans $B_n(2r)$, \'egale \`a 1 sur $B_n(r)$ et
telle que $\ddc\chi\leq c''r^{-2}\omega$ avec $c''>0$ o\`u on a pos\'e
$\omega:=\ddc|z|^2$. On a 
\begin{eqnarray*}
\mu(B_n(r)) & \leq & \int \chi_n \ddc G_m \wedge T_{k-1} = \int_{B_n(2r)} \ddc
\chi_n G_m\wedge T_{k-1} \\
& \leq & \int_{B_n(2r)} c''r^{-2} \omega c r^{a^*_m} c''\wedge
T_{k-1} = c_1 r^{a^*_m-2} \int_{B_n(2r)} \omega\wedge T_{k-1}
\end{eqnarray*}
o\`u $c_1:=cc''$. 
\par
On peut r\'ep\'eter ce proc\'ed\'e $(k-1)$ fois et on obtient, pour
des constantes $c_k>0$ et $c_k'>0$ convenables, que
$$\mu(B_n(r))\leq c_k r^{a^*-2k}\int_{B(2^kr)} \omega^k =c_k' r^{a^*}.$$
Par cons\'equent, 
$$\mu(\Sigma)\leq \sum \mu(B_n(r))\leq \epsilon r^{-a}c_k' r^{a^*} =
c_k'\epsilon.$$ 
Ceci est vrai pour tout $\epsilon>0$. Donc $\mu(\Sigma)=0$.
\end{preuve}
\begin{remarque} \rm
Les r\'esultats pour les endomorphismes $(\pi,s)$-r\'eguliers
s'appliquent  aux automorphismes, bien s\^ur on a alors $s\leq m$. 
Si $f$ est $(\pi,s)$-r\'egulier et
$f^{-1}$ est $(\tilde\pi,\tilde s)$-r\'egulier on construit des courants
invariants $T^+$ pour $f$ et $T^-$ pour $f^{-1}$. On peut, pour $\pi$,
$\tilde\pi$, $s$ et $\tilde s$ convenables, consid\'erer
la mesure invariante $\mu:=T^+\wedge T^-$. Le cas le plus simple de
cette situation est celui des applications de H\'enon. On trouve
d'autres exemples dans \cite{Sibony2} et \cite{GuedjSibony}.
\end{remarque}
\section{D'autres remarques}
Soit $f:\C^k\longrightarrow\C^k$ une application polynomiale propre
de degr\'e topologique $d_t\geq 2$. Supposons que l'infini soit
attirant dans le sens o\`u il existe $c>1$ tel que 
$|f(z)|\geq c|z|$ pour $|z|$ grand. On peut construire la mesure d'\'equilibre
$\mu$ de $f$ comme la limite faible de la suite $d_t^{-n}(f^n)^*\nu$ o\`u
$\nu$ est une mesure de probabilit\'e qui ne charge pas les ensembles
pluripolaires \cite{DinhSibony2}. 
La mesure $\mu$ est m\'elangeante et ne d\'epend pas de $\nu$.
Lorsque l'exposant de Lojasiewicz $\lambda_1$ de $f$ est strictement
sup\'erieur \`a 1, en utilisant la m\'ethode de Lyubich \cite{Lyubich} et
de Briend-Duval \cite{BriendDuval2, BriendDuval1}, on
montre \cite{DinhSibony2}
 que $d_t^{-n}(f^n)^*\delta_z$ tend vers $\mu$ pour tout $z$
hors d'un ensemble exceptionnel $\E$ qui est analytique. On a not\'e
$\delta_z$ la masse de Dirac en $z$. Les exposants de Lyapounov de
$\mu$ sont minor\'es par $\log\lambda_1/2$. Les points p\'eriodiques
r\'epulsifs sont denses et \'equidistribu\'es sur $\supp(\mu)$. La
vitesse de m\'elange de $\mu$ est de l'ordre $\lambda_1^{-n}$. Cette
mesure $\mu$ est de plus l'unique mesure d'entropie maximale $\log
d_t$. Pour cette derni\`ere propri\'et\'e de $\mu$, 
il suffit de reprendre la preuve de Lyubich \cite{Lyubich} et Briend-Duval
\cite{BriendDuval2} en rempla\c cant un calcul cohomologique par le
lemme de comparaison suivant:
\begin{lemme} Soit $\Omega$ une forme de bidegr\'e $(k-1,k-1)$
positive ferm\'ee dans $\C^k$. 
Alors on a pour tout $m\geq 0$
$$\int\Omega\wedge (f^m)^*\omega\leq\frac{1}{\lambda_1}\int\Omega\wedge
(f^{m+1})^*\omega.$$
\end{lemme}
\begin{preuve} Il suffit d'appliquer la proposition 5.4 pour
$V=\C^k$, $\rho(z)=\log(1+|z|^2)$, 
$v_1=\log(1+|f^m|^2)-A$, $v_2=\lambda_1^{-1}\log(1+|f^{m+1}|^2)$ et $A$
une constante suffisamment grande.
\end{preuve}
\par
Les r\'esultats cit\'es ci-dessus sont valables, en particulier,
si $f$ est
semi-r\'egulier avec $\alpha_m>1$. Dans ce cas, la mesure
d'\'equilibre $\mu$ est \'egale \`a $T_k$ qui est une intersection
g\'en\'eralis\'ee de courants positifs ferm\'es. En particulier,
d'apr\`es la proposition 5.3, 
les fonctions p.s.h. sont $\mu$-int\'egrables: c'est une
mesure PLB \cite{DinhSibony2}. 
\par
En g\'en\'eral, une application $f$ d'exposant de Lojasiewicz
$\lambda_1>1$ n'est pas conjugu\'ee \`a une application
semi-r\'eguli\`ere. Il peut exister une infinit\'e de vitesses
d'\'echapement vers l'infini et on rencontre 
d'autres ph\'enom\`enes dynamiques. 
Ceci fera l'objet d'un prochain travail avec R. Dujardin.
\par
Dans la suite, nous consid\'erons des exemples d'applications
r\'eguli\`eres avec $d_m=1$. 
Soit $f:\C^2\longrightarrow \C^2$ d\'efini par
$f(z)=(P(z),a z_1+bz_2)$ o\`u $a$, $b$ sont des nombres complexes,
$|b|>1$, $P$ est un polyn\^ome de degr\'e $d\geq 2$. Supposons que le
coefficient de $z_1^d$ dans $P$ soit non nul. L'application $f$ 
est r\'eguli\`ere et de
degr\'e topologique $d_t=d$. En particulier, 
elle est alg\'ebriquement stable.
D'apr\`es le th\'eor\`eme 2.1, l'exposant de Lojasiewicz et l'exposant
de Lojasiewicz asymptotique de $f$
sont \'egaux \`a 1.
On v\'erifie $|f(z)|\geq c |z|$
pour tout $1<c<|b|$ fix\'e et pour $|z|$ assez
grand.
\par 
Nous ne savons pas si $\E$ est analytique et si
les exposants de Lyapounov sont strictement positifs lorsque
$P(0,z_2)\not\equiv 0$.
\par 
On peut construire la fonction de Green p.s.h., continue, positive,
invariante: $G_1\circ f=dG_1$ et le courant de Green $T_1:=\ddc G_1$, 
positif, ferm\'e, invariant par $f$:
$f^*T_1=dT_1$.
L'ensemble $\K_1:=\{G_1=0\}$, qui est le compl\'ement du bassin
d'attraction $U_1$ de
$X_1:=[1:0:0]$, n'est pas compact et donc
n'est pas \'egal \`a $\K$. Le seul point adh\'erent
\`a $\K_1$ dans l'hyperplan \`a l'infini est l'unique point
d'ind\'etermination $I_1:=[0:1:0]$. On en d\'eduit que 
si $z\in \K_1$ tend vers $I_1$, on a $|z_1|\leq A|z_2|^{(d-1)/d}$
car $f(\K_1)=\K_1$. On a aussi  $|f(z)|\simeq
|b||z|$. 
Le courant $T_1$ ne charge pas les ensembles
pluripolaires car son potentiel $G_1$ est continu. Notons que
$T_1\wedge T_1=0$ dans $\C^2$ car $f^*(T_1\wedge T_1)=d^2T_1\wedge
T_1$ et $d_t<d^2$. 
On a, d'apr\`es \cite[th\'eor\`eme 3.2.1]{DinhSibony2}  
$$\mu=\lim (f^n)^*\omega\wedge T_1 = \lim \ddc u_n\wedge T_1$$
o\`u $u_n(z):=\log^+|f^n(z)| -n\log |b|$. 
\par
Il est clair que $\lim u_n=-\infty$ sur $\K$. La suite $(u_n)$
converge uniform\'ement sur les compacts de $\K_1\setminus \K$. En
effet, on a, en posant $f^n(z)=w=(w_1,w_2)$
\begin{eqnarray*}
|u_{n+1}(z)-u_n(z)| & \sim & \log \left|\frac{a w_1+bw_2}{bw_2}\right|
  \sim \log \left|1+\frac{a w_1}{bw_2}\right|\\
& \sim & \left|\frac{w_1}{w_2}\right|
\leq \frac{1}{|w|^{1/d}} \leq \frac{A}{|c|^{n/d}}.
\end{eqnarray*}
La constante $c$ \'etant sup\'erieure \`a 1, 
la fonction $u:=\lim u_n$ est donc continue sur
  $\K_1\setminus \K$ et v\'erifie la relation $u\circ f=u+\log |b|$. 
Elle est $T_1$-p.s.h. sur $\K_1\setminus K$ et tend vers $-\infty$
quand $z$ tend vers $\K$.
\section{Appendice: fonctions T-p.s.h.}
Nous explicitons dans cet appendice 
quelques propri\'et\'es, que nous utilisons, des fonctions
$T$-p.s.h., \cad 
les fonctions p.s.h. relativement \`a un courant positif ferm\'e $T$. 
Ces propri\'et\'es sont
classiques dans le cadre des fonctions p.s.h.
et les d\'emonstrations sont de simples extensions du cas des
fonctions p.s.h.
\cite{BedfordTaylor}, \cite{Taylor}, 
\cite{Demailly}, \cite{FornaessSibony3}.  
\par
Soit $T$ un courant positif ferm\'e de bidegr\'e $(j,j)$ dans une
vari\'et\'e k\"ahl\'erienne $V$ de dimension $k\geq 1$ (\voir
\cite{Demailly}, \cite{GriffithsHarris}, \cite{Lelong} 
pour les d\'efinitions de base). Rappelons
cependant que pour un courant 
$T\geq 0$ de bidegr\'e $(1,1)$, on a
localement $T=\ddc u$ o\`u $u$ est une fonction p.s.h. On dit que $u$
est un {\it potentiel local} de $T$.
Rappelons les d\'efinitions de \cite{BerndtssonSibony}.
Une fonction semi-continue
sup\'erieurement (s.c.s) $v$ sur $\supp(T)$ est dite {\it $T$-p.s.h.} 
si elle
est localement limite d\'ecroissante d'une suite $(v^{(n)})$ de fonctions
${\cal C}^2$ v\'erifiant $\ddc v^{(n)}\wedge T\geq 0$. On dira que $v$ est
{\it fortement $T$-p.s.h.}
si les $v^{(n)}$ sont p.s.h. au voisinage du point
consid\'er\'e de
$\supp(T)$.
On pose $$\sigma_T:=\frac{T\wedge \omega^{k-j}}{(k-j)!}$$ o\`u
$\omega$ d\'esigne la forme de K\"ahler sur $V$. C'est {\it la mesure trace
de $T$}. On pose $\|T\|_K:=\sigma_T(K)$.
%
%
%
\begin{proposition} Soient $L_1\subset\subset L_2$  deux
  compacts de $V$. Soient 
$v\in\Loneloc(\sigma_T)$ une fonction $T$-p.s.h. et 
$v_1$, $\ldots$, $v_q$ des fonctions $T$-p.s.h., localement born\'ees.
Alors
\begin{enumerate}
\item Le courant
  $\ddc v\wedge T:=\ddc(vT)$ est positif
ferm\'e de bidegr\'e $(j+1,j+1)$. 
\item
Il existe une constante $c_{L_1,L_2}>0$, ind\'ependante de $v$
et des $v_i$, telle que 
$\|\ddc v\wedge T\|_{L_1}\leq c_{L_1,L_2} \|vT\|_{L_2}$
et
$$\|\ddc v_1\wedge \ldots \wedge \ddc v_q\wedge T\|_{L_1} \leq
c_{L_1,L_2} \|v_1\|_{\Linfty(L_2)}\ldots \|v_q\|_{\Linfty(L_2)}
\|T\|_{L_2}.$$
\end{enumerate}
\end{proposition}
\begin{preuve} Pour la positivit\'e, il suffit d'observer que localement 
$$\ddc (vT)=\lim \ddc(v^{(n)}T) =\lim \ddc v^{(n)}\wedge T\geq 0.$$
\par
Soit  $\chi$ une fonction de classe ${\cal
  C}^2$ \`a support dans $L_2$ et \'egale \`a 1 sur $L_1$.
Posons $c_{L_1,L_2}:= \|\ddc \chi\|_{\Linfty(L_2)}$.
Il suffit d'estimer par int\'egration par
parties:
$$\|\ddc v\wedge T\|_{L_1} \leq \|\chi \ddc(vT)\|_{L_2}
\leq \|\ddc\chi vT\|_{L_2} \leq c_{L_1,L_2}
\|vT\|_{L_2}.$$
\par
La derni\`ere relation se d\'emontre par r\'ecurrence sur $q$.
\end{preuve}
\par
Le r\'esultat suivant 
est l'analogue du th\'eor\`eme de continuit\'e
de Bedford-Taylor \cite{BedfordTaylor} (\voir aussi \cite{Demailly}).
\begin{proposition} Soient 
  $v_1$, $\ldots$, $v_q$ des fonctions $T$-p.s.h. localement born\'ees.
Soient $v_1^{(n)}$, $\ldots$, $v_q^{(n)}$ des fonctions $T$-p.s.h.
  d\'ecroissant vers $v_1$, $\ldots$, $v_q$. Alors 
\begin{enumerate}
\item $v_1^{(n)}\ddc v_2^{(n)}\wedge\ldots\wedge \ddc v_q^{(n)} \wedge
  T\rightharpoonup  v_1\ddc v_2\wedge\ldots\wedge \ddc v_q \wedge
  T$ faiblement.
\item $\ddc v_1^{(n)}\wedge\ddc v_2^{(n)}\wedge\ldots\wedge \ddc
  v_q^{(n)} \wedge
  T\rightharpoonup  \ddc v_1\wedge \ddc v_2\wedge\ldots\wedge \ddc v_q \wedge
  T$ faiblement.
\end{enumerate}
En particulier, l'application
$(v_1,\ldots,v_q)\mapsto \ddc v_1\wedge\ldots\wedge \ddc v_q\wedge T$
est sym\'etrique en $v_1$, $\ldots$, $v_q$.
Si l'on suppose que les $v_l$ et $v_l^{(n)}$ sont fortement $T$-p.s.h., la
d\'ecroissance est superflue, il suffit de supposer $v_l^{(n)}\geq v_l$
pour $1\leq l\leq q$.
\end{proposition}
\begin{preuve} Pour reprendre la d\'emonstration de
  \cite{BedfordTaylor} et
  \cite{Demailly}, il suffit
  de faire les remarques suivantes.
\par
Soit $\chi$ une fonction convexe dans $\R^q$ croissante par
rapport \`a chaque variable. Si $v_1$, $\ldots$, $v_q$ sont $T$-p.s.h.
alors
$\chi(v_1,\ldots,v_q)$ est $T$-p.s.h. 
En particulier, le sup d'un
nombre fini de fonctions $T$-p.s.h. 
l'est aussi. Cela permet de se
ramener au cas o\`u $V$ est la boule unit\'e et 
toutes les fonctions $v^{(n)}_l$ sont \'egales \`a
$|z|^2-1$ au voisinage de la sph\`ere unit\'e. On peut alors appliquer
les int\'egrations par parties usuelles comme dans
\cite{BedfordTaylor} ou \cite{Demailly}.
\par
Lorsque les fonctions sont fortement $T$-p.s.h., une utilisation du lemme
de Hartogs comme dans \cite[p.405]{FornaessSibony3} permet de montrer le
r\'esultat. 
\end{preuve}
\begin{proposition}{\bf (Chern-Levine-Nirenberg)}
Soient $L_1\subset\subset L_2$ deux compacts de
  $V$.
Soient $v_1$, $\ldots$, $v_q$ des fonctions $T$-p.s.h. 
localement born\'ees.
Alors il existe une constante $c_{L_1,L_2}>0$
  telle que pour toute fonction $T$-p.s.h. 
$\varphi\in\Loneloc(\sigma_T)$, on ait
$$\|\varphi \ddc v_1\wedge \ldots \wedge \ddc v_q \wedge T\|_{L_1}\leq
  c_{L_1,L_2} \|\varphi\|_{\Lone(\sigma_T,L_2)} \|v_1\|_{\Linfty(L_2)}\ldots
  \|v_q\|_{\Linfty(L_2)}.$$
En particulier, si $T=\ddc u_j\wedge\ldots\wedge \ddc u_1$ avec $u_1$
  p.s.h. born\'ee et $u_l$
  fonction $(\ddc u_{l-1}\wedge\ldots\wedge \ddc u_1)$-p.s.h. localement
  born\'ee pour tout $2\leq l\leq j$, alors les fonctions p.s.h. sont
  localement $\sigma_T$-int\'egrables et $T$ ne charge pas les ensembles
  pluripolaires. 
\end{proposition}
\begin{preuve}
Il suffit de reprendre la
d\'emonstration de \cite[p.126]{Demailly} en introduisant le courant
$T$ (\voir \'egalement \cite{FavreGuedj}). Pour la commodit\'e du
lecteur, nous donnons ici la preuve.
On se ram\`ene au cas o\`u $L_1$ et $L_2$ sont des boules
centr\'ees en $0$ et de rayon respectif $R'$ et $R$ 
et o\`u tous les $v_j$ sont \'egaux \`a $|z|^2-R^2$ pour $R_1<|z|<R$
avec un $R'<R_1<R$. On peut aussi supposer que $j+q=k-1$. Soit $0\leq
\chi\leq R^2$ une fonction \'egale \`a $R^2-|z|^2$ pour $|z|<R'$ et \`a
support dans $B_{R_1}$ de centre $0$ et de rayon $R_1$. 
On peut supposer $\varphi<0$. On a pour un $c>0$ 
ind\'ependant de $\varphi$
\begin{eqnarray*}
I & := & \int_{|z|<R'} 
-\varphi \ddc v_1\wedge \ldots \wedge \ddc v_q \wedge T
\wedge \ddc |z|^2\\
& = & \int_{|z|<R_1} \varphi \ddc v_1\wedge \ldots
\wedge \ddc v_q \wedge T \wedge \ddc\chi - \\
& & - \int_{R'<|z|<R_1} 
\varphi \ddc v_1\wedge \ldots
\wedge \ddc v_q \wedge T \wedge \ddc\chi \\
&\leq & 
\int \chi \ddc\varphi \wedge  \ddc v_1\wedge \ldots
\wedge \ddc v_q \wedge T + \\
& & 
+ c \|\varphi\|_{\Lone(\sigma_T,L_2)}
\|v_1\|_{\Linfty(L_2)}\ldots
  \|v_q\|_{\Linfty(L_2)}
\end{eqnarray*}
car $\ddc v_j=\ddc |z|^2$ sur le support de $\ddc\chi$. Par suite,
\begin{eqnarray*}
I &\leq & 
R^2\int_{|z|<R_1} \ddc\varphi \wedge  \ddc v_1\wedge \ldots
\wedge \ddc v_q \wedge T + \\
& & 
+ c \|\varphi\|_{\Lone(\sigma_T,L_2)}
\|v_1\|_{\Linfty(L_2)}\ldots
  \|v_q\|_{\Linfty(L_2)}.
\end{eqnarray*}
Il reste \`a majorer la derni\`ere int\'egrale. Pour ceci, puisque
$\varphi$ est la limite d\'ecroissante de fonctions $T$-p.s.h. lisses, on
peut supposer $\varphi$ lisse et passer ensuite \`a la limite. 
Soit $R_2$ tel que $R_1<R_2<R$. 
D'apr\`es la proposition 5.2, 
on a 
\begin{eqnarray*}
\lefteqn{\int_{|z|<R_1} \ddc\varphi \wedge  \ddc v_1\wedge \ldots
\wedge \ddc v_q \wedge T} & &\\
&\hspace{2cm}  = & \int_{|z|<R_1} \ddc v_1\wedge \ldots
\wedge \ddc v_q \wedge \ddc (\varphi T) \\
&\hspace{2cm}  \leq & c_1\|v_1\|_{\Linfty(L_2)}\ldots
  \|v_q\|_{\Linfty(L_2)} \|\ddc \varphi T\|_{\Lone(B_{R_2})}\\
&\hspace{2cm}  \leq & c_2\|v_1\|_{\Linfty(L_2)}\ldots
  \|v_q\|_{\Linfty(L_2)} \|\varphi\|_{\Lone(\sigma_T,L_2)}.
\end{eqnarray*}
Les  derni\`eres in\'egalit\'es sont des cons\'equences de la proposition
5.1; $c_1$, $c_2$ sont des constantes ind\'ependantes de $\varphi$.
\end{preuve}
\par
Enon\c cons un {\it th\'eor\`eme de comparaison} 
\index{th\'eor\`eme!de comparaison}
(\voir \cite{Taylor}).
\begin{proposition} Soit $V$ une vari\'et\'e de Stein de dimension
  $k\geq 1$.
Soit $\rho$ une fonction p.s.h. d'exhaustion de $V$, \ie $\rho$
  tend vers l'infini \`a l'infini. Posons $\omega:=\ddc \rho$. 
Soit $T$ un courant de bidegr\'e $(k-1,k-1)$ positif ferm\'e. Soient
$v_1$, $v_2$ deux fonctions $T$-p.s.h. avec $v_j\in \Loneloc(\sigma_T)$. 
On suppose que l'une des deux conditions suivantes est v\'erifi\'ee
\begin{enumerate}
\item Sur $\supp(T)$, on a 
$v_1\leq v_2$ \`a l'infini et $v_2\rightarrow\infty$ \`a l'infini.
\item Sur $\supp(T)$, on a $v_1\leq v_2$ \`a l'infini, 
$\int T\wedge \ddc\rho<\infty$ et $v_2+\epsilon\rho\rightarrow \infty$ \`a
l'infini pour tout $\epsilon>0$. 
\end{enumerate}
Alors
$$\int_V\ddc v_1\wedge T\leq \int_V\ddc v_2 \wedge T.$$
\end{proposition}
\begin{preuve} D'apr\`es la proposition 5.2,
il suffit de montrer la 
proposition pour $v_{i,M}:=\max(v_i,-M)$ puis de faire
tendre $M$ vers $+\infty$. On peut donc supposer que les $v_i$ sont
localement born\'ees. Soient $\epsilon>0$ assez petit 
et $R>0$ assez grand. Posons
pour le premier cas
$$W_\epsilon:=\max \big(v_1+A,(1+\epsilon)v_2\big)$$
et pour le second cas
$$W_\epsilon:=\max(v_1+A,v_2+\epsilon\rho).$$
La constante $A$ est telle que sur $(v_2<R)$
(resp. sur $(v_2+\epsilon\rho<R)$ pour le second cas) 
on ait $W_\epsilon=v_1+A$. Soit $\chi$ une fonction test \`a support compact,
$0\leq \chi\leq 1$, $\chi=1$ sur $(v_2< R)$
(resp. sur $(v_2+\epsilon\rho<R)$) et telle que $(\d\chi\not =0)$
soit contenu dans l'ensemble o\`u $W_\epsilon=(1+\epsilon)v_2$ (resp. 
$W_\epsilon:= v_2+\epsilon \rho$). On a pour
le premier cas
\begin{eqnarray*}
\int_{(v_2 <R)} \ddc v_1\wedge T & = & \int_{(v_2<R)} \ddc
W_\epsilon \wedge T \leq \int \chi \ddc W_\epsilon \wedge T
\\
&= & \int W_\epsilon
\ddc \chi \wedge T  =  \int (1+\epsilon)v_2 \ddc \chi \wedge T \\
& =  & (1+\epsilon)
\int \ddc\chi \wedge v_2 T = (1+\epsilon)\int \chi \ddc v_2\wedge T\\
& \leq & (1+\epsilon) \int_V\ddc v_2\wedge T.
\end{eqnarray*}
On peut faire tendre $\epsilon$ vers $0$ puis $R$ vers l'infini. Le
second cas se d\'emontre de la m\^eme mani\`ere.
\end{preuve}
\small

\par\noindent
Tien-Cuong Dinh et Nessim Sibony\\
Math\'ematique - B\^at. 425, UMR 8628\\
Universit\'e Paris-Sud, 91405 Orsay, France.\\ 
E-mails: Tiencuong.Dinh@math.u-psud.fr, Nessim.Sibony@math.u-psud.fr
\end{document}